\documentclass[11pt,a4paper]{article}
\usepackage{mathrsfs}
\usepackage{amsmath}
\usepackage{amsfonts}
\usepackage{amssymb}

\usepackage{graphicx,subfigure,epsfig}
\usepackage{enumerate}
\usepackage{arydshln}
\usepackage[pdftex,linkcolor=blue,citecolor=blue,backref=pagebackref]{hyperref}
\usepackage{color}
\usepackage{pdfpages}
\usepackage{subfigure}
\usepackage{tabularx}
\allowdisplaybreaks[4]

\begin{document}

\renewcommand{\thefootnote}{\fnsymbol{footnote}}

\newtheorem{theorem}{Theorem}[section]
\newtheorem{corollary}[theorem]{Corollary}
\newtheorem{definition}[theorem]{Definition}
\newtheorem{conjecture}[theorem]{Conjecture}
\newtheorem{question}[theorem]{Question}
\newtheorem{lemma}[theorem]{Lemma}
\newtheorem{proposition}[theorem]{Proposition}
\newtheorem{example}[theorem]{Example}
\newtheorem{fact}[theorem]{Fact}
\newenvironment{proof}{\noindent {\bf
Proof.}}{\rule{3mm}{3mm}\par\medskip}
\newcommand{\remark}{\medskip\par\noindent {\bf Remark.~~}}
\newcommand{\pp}{{\it p.}}
\newcommand{\de}{\em}
\newcommand{\qbinom}[2]{\genfrac{[}{]}{0pt}{}{#1}{#2}}

\newcommand{\JEC}{{\it Europ. J. Combinatorics},  }
\newcommand{\JCTB}{{\it J. Combin. Theory Ser. B.}, }
\newcommand{\JCT}{{\it J. Combin. Theory}, }
\newcommand{\JGT}{{\it J. Graph Theory}, }
\newcommand{\ComHung}{{\it Combinatorica}, }
\newcommand{\DM}{{\it Discrete Math.}, }
\newcommand{\ARS}{{\it Ars Combin.}, }
\newcommand{\SIAMDM}{{\it SIAM J. Discrete Math.}, }
\newcommand{\SIAMADM}{{\it SIAM J. Algebraic Discrete Methods}, }
\newcommand{\SIAMC}{{\it SIAM J. Comput.}, }
\newcommand{\ConAMS}{{\it Contemp. Math. AMS}, }
\newcommand{\TransAMS}{{\it Trans. Amer. Math. Soc.}, }
\newcommand{\AnDM}{{\it Ann. Discrete Math.}, }
\newcommand{\NBS}{{\it J. Res. Nat. Bur. Standards} {\rm B}, }
\newcommand{\ConNum}{{\it Congr. Numer.}, }
\newcommand{\CJM}{{\it Canad. J. Math.}, }
\newcommand{\JLMS}{{\it J. London Math. Soc.}, }
\newcommand{\PLMS}{{\it Proc. London Math. Soc.}, }
\newcommand{\PAMS}{{\it Proc. Amer. Math. Soc.}, }
\newcommand{\JCMCC}{{\it J. Combin. Math. Combin. Comput.}, }
\newcommand{\GC}{{\it Graphs Combin.}, }
\newcommand{\LAA}{{\it Linear Algeb. Appli.}, }

\title{ \bf $\mathcal{L}$-intersecting or Configuration Forbidden Families
\\on Set Systems and Vector Spaces over Finite Fields}

\author{Jiuqiang Liu$^{a,b,} \thanks{The corresponding author}$, Guihai Yu$^{a, *}$, Lihua Feng$^{c}$, and Yongtao Li$^{c}$\\
{\small  $^a$ College of Big Data Statistics, Guizhou University of Finance and Economics}\\
 {\small Guiyang, Guizhou, 550025, China}\\
 {\small $^{b}$ Department of Mathematics, Eastern Michigan University}\\
{\small Ypsilanti, MI 48197, USA}\\
{\small $^{c}$ School of Mathematics and Statistics, HNP-LAMA, Central South University}\\
{\small Changsha, Hunan, 410083, China}\\
{\small E-mail: { \tt jiuqiang68@126.com, yuguihai@mail.gufe.edu.cn, fenglh@163.com}}}

\maketitle
\vspace{-0.5cm}

\begin{abstract}
In this paper, we derive a tight upper bound for the size of an intersecting $k$-Sperner family of subspaces of the $n$-dimensional vector space $\mathbb{F}_{q}^{n}$ over finite field $\mathbb{F}_{q}$ which gives a $q$-analogue of the Erd\H{o}s' $k$-Sperner Theorem, and we then establish a general relationship between upper bounds for the sizes of families of subsets of $[n] = \{1, 2, \dots, n\}$ with property $P$ and upper bounds for the sizes of families of subspaces of $\mathbb{F}_{q}^{n}$ with property $P$, where $P$ is either $\mathcal{L}$-intersecting or forbidding certain configuration. Applying this relationship, we derive generalizations of the well known results about the famous Erd\H{o}s matching conjecture and Erd\H{o}s-Chv\'atal simplex conjecture to linear lattices.

\end{abstract}

\noindent
{{\bf Key words:}  Erd\H{o}s-Ko-Rado theorem, Erd\H{o}s matching conjecture;  Erd\H{o}s-Chv\'atal simplex conjecture; hereditary property; $\mathcal{L}$-intersecting family; Sperner family} \\
{{\bf AMS Classifications:} 05D05. } \vskip 0.1cm

\section{Introduction}
\hspace*{0.5cm}
Throughout the paper, we denote  $[n]=\{1, 2, \dots, n\}$. A family $\mathcal{F}$ of subsets of $[n]$ is called
$t$-$intersecting$ if $|E \cap F| \geq t$ for every pair of distinct subsets $E, F \in
\mathcal{F}$, and it is called {\em intersecting} when $t = 1$. Let $\mathcal{L}=\{l_1, l_2, \dots, l_s\}$ be a set of $s$ nonnegative integers
with $0 \leq l_{1} < \cdots < l_{s}$.
A family $\mathcal{F}$ of subsets of $[n]$ is called
$\mathcal{L}$-$intersecting$ if $|E \cap F| \in \mathcal{L}$ for
every pair of subsets $E$, $F$ in $\mathcal{F}$.
$\mathcal{F}$ is $k$-{\em uniform} if it is a collection of $k$-subsets
of $[n]$. Thus, a $k$-uniform intersecting family is
$\mathcal{L}$-intersecting for $\mathcal{L}=\{1, 2, \dots, k-1\}$.
A $k$-uniform $t$-intersecting family is a special $\mathcal{L}$-intersecting family with $\mathcal{L}=\{t, t+1, \dots, k-1\}$.
A {\em Boolean lattice} $\mathcal{B}_{n}$ is the set $2^{[n]}$ of all subsets of $[n]$ with the ordering being the containment relation and a {\em linear lattice}
$\mathcal{L}_{n}(q)$ is the set of all subspaces of the $n$-dimensional vector space $\mathbb{F}_{q}^{n}$ over the finite field $\mathbb{F}_{q}$ with the ordering being the relation of inclusion of subspaces.

In 1961, Erd\H{o}s, Ko, and Rado \cite{ekr} proved the following classical result.

\begin{theorem}\label{thm1.1} (Erd\H{o}s, Ko, and Rado, \cite{ekr}).
Let $n \geq 2k$, and let $\mathcal{A}$ be a $k$-uniform intersecting
family of subsets of $[n]$. Then $|\mathcal{A}|\leq {{n-1} \choose
{k-1}}$, with equality only when $\mathcal{F}$ consists of all
$k$-subsets containing a common element.
\end{theorem}

Since then, many intersection theorems have appeared in the literature, see, for example, the next three well-known theorems.

\begin{theorem}\label{thm1.2} (Ray-Chaudhuri and Wilson, \cite{rw}).
Let $\mathcal{L}=\{l_1, l_2, \dots, l_s\}$ be a set of $s$
nonnegative integers. If $\mathcal{A}$ is a $k$-uniform
$\mathcal{L}$-intersecting family of subsets of $[n]$, then
$|\mathcal{A}|\leq {{n} \choose {s}}$.
\end{theorem}

\begin{theorem}\label{thm1.3} (Frankl and Wilson, \cite{fw}). Let
$\mathcal{L}=\{l_1, l_2, \dots, l_s\}$ be a set of $s$ nonnegative
integers. If $\mathcal{A}$ is an $\mathcal{L}$-intersecting family
of subsets of $[n]$, then
\[|\mathcal{A}|\leq {{n} \choose {s}}+ {{n}
\choose {s-1}}+ \cdots + {{n} \choose {0}}.\]
\end{theorem}

In terms of the parameters $n$ and $s$, the inequality in Theorem \ref{thm1.3} is best possible, as shown by the set of all subsets of $[n]$ with sizes at most $s$ and taking $\mathcal{L}=\{0, 1, \dots, s - 1\}$.

\begin{theorem}\label{thm1.4} (Alon, Babai, and Suzuki, \cite{abs}). Let
$\mathcal{L}=\{l_1, l_2, \dots, l_s\}$ be a set of $s$ nonnegative
integers and $K = \{k_1, k_2, \dots, k_r \}$ be a set of integers
satisfying $k_i > s-r$ for every $i$.
Suppose that $\mathcal{A} = \{A_1, A_2, \dots, A_{m}\}$ is a family of subsets of $[n]$ such that $|A_i| \in K$ for every $1 \leq i \leq m$
and $|A_{i} \cap A_{j}| \in \mathcal{L}$ for every pair $i \neq j$.
Then
\[m \leq {{n} \choose {s}}+ {{n} \choose {s-1}}+ \cdots + {{n} \choose {s-r+1}}.\]
\end{theorem}

Recall that a family of subsets of $[n]$ is called a {\em Sperner family} (or antichain) if there
are no two different members of the family such that one of them contains the
other, and a family $\mathcal{F}$ of subsets of $[n]$ is $k$-$Sperner$
if all chains in $\mathcal{F}$ have length at most $k$. Define $\sum (n, k)$ to be the sum of the
$k$ largest binomial coefficients of order $n$, i.e., $\sum (n, k) = \sum_{i = 1}^{k}{{n} \choose {\lfloor \frac{n - k}{2}\rfloor + i}}$.
Let $\sum^{*} (n, k)$ be the collection of families consisting of the corresponding full levels,
i.e., if $n + k$ is odd, then $\sum^{*} (n, k)$ contains one family $\cup_{i = 1}^{k}{{[n]} \choose {\lfloor \frac{n - k}{2}\rfloor + i}}$
(where ${{[n]} \choose {k}}$ is the set of all $k$-subsets of $[n]$);
if $n + k$ is even, then $\sum^{*} (n, k)$ contains two families of the same size $\cup_{i = 0}^{k - 1}{{[n]} \choose {\frac{n - k}{2}+ i}}$ and
$\cup_{i = 1}^{k}{{[n]} \choose {\frac{n - k}{2}+ i}}$. The following theorem by Erd\H{o}s \cite{e} generalizes the classical Sperner theorem.

\begin{theorem}\label{thm1.5} (Erd\H{o}s, \cite{e}). Suppose that $\mathcal{A}$ is a $k$-Sperner family
of subsets of $[n]$. Then
\[|\mathcal{A}|\leq \sum (n, k).\]
Moreover, if $|\mathcal{A}| = \sum (n, k)$, then $\mathcal{A} \in \sum^{*} (n, k)$.
\end{theorem}

Throughout this paper, we denote the Gaussian (or $q$-binomial) coefficient of order $n$ by
\begin{equation}
\qbinom{n}{k} = \prod_{0 \leq i \leq k - 1}\frac{q^{n - i} - 1}{q^{k - i} - 1}.
\label{1.1}\end{equation}
It is well known that the number of all $k$-dimensional
subspaces of $\mathbb{F}_{q}^{n}$ is equal to $\left[n \atop k \right]$.

Results on the sizes of families of subspaces of $\mathbb{F}_{q}^{n}$ analogous to those theorems on $\mathcal{L}$-intersecting families of subsets of $[n]$ have appeared in the literature. The next generalization of  Erd\H{o}s-Ko-Rado Theorem (Theorem \ref{thm1.1}) is proved by Greene and Kleitman \cite{gk} and Hsieh \cite{h}.

\begin{theorem}\label{thm1.6} (Greene and Kleitman \cite{gk}, Hsieh \cite{h}).
Let $n \geq 2k + 1$. Suppose that $\mathcal{V}$ is a family of $k$-dimensional subspaces of $\mathbb{F}_{q}^{n}$ satisfying that $\dim(V_{i} \cap V_{j}) > 0$ for any distinct subspaces $V_{i}$ and $V_{j}$ in $\mathcal{V}$.
Then
\[|\mathcal{V}| \leq \qbinom{n-1}{k-1}\]
and equality holds if and only if there exists a $1$-dimensional subspace $R$ of $\mathbb{F}_{q}^{n}$ such that $R \subseteq V$ for every $V \in \mathcal{V}$.
\end{theorem}

Frankl and Graham \cite{fg} derived the following theorem on $\mathcal{L}$-intersecting family of $k$-dimensional subspaces of $\mathbb{F}_{q}^{n}$.

\begin{theorem}\label{thm1.7} (Frankl and Graham \cite{fg}).
Let $\mathcal{L}=\{l_1, l_2, \dots, l_s\}$ be a set of $s$ nonnegative integers. Suppose that $\mathcal{V}$ is a family of $k$-dimensional subspaces of $\mathbb{F}_{q}^{n}$ such that $\dim(V_{i} \cap V_{j}) \in \mathcal{L}$ for any distinct subspaces $V_{i}$ and $V_{j}$ in $\mathcal{V}$.
Then
\[|\mathcal{V}| \leq \qbinom{n}{s}.\]
\end{theorem}

In 1990, Lefmann \cite{l} proved the following $\mathcal{L}$-intersecting theorem for ranked finite lattices,
and in 2001, Qian and Ray-Chaudhuri \cite{qr} extended to quasi-polynomial semi-lattices.

\begin{theorem}\label{thm1.8} (Lefmann \cite{l}, Qian and Ray-Chaudhuri \cite{qr}).
Let $\mathcal{L}=\{l_1, l_2, \dots, l_s\}$ be a set of $s$ nonnegative integers. Suppose that $\mathcal{V}$ is a family of subspaces of $\mathbb{F}_{q}^{n}$ satisfying that $\dim(V_{i} \cap V_{j}) \in \mathcal{L}$ for any distinct subspaces $V_{i}$ and $V_{j}$ in $\mathcal{V}$.
Then
\[|\mathcal{V}| \leq \qbinom{n}{s}
+ \qbinom{n}{s-1} + \cdots + \qbinom{n}{0}.\]
\end{theorem}

Alon et al. \cite{abs} derived the next result in 1991.

\begin{theorem}\label{thm1.9} (Alon, Babai, and Suzuki \cite{abs}).
Let $\mathcal{L}=\{l_1, l_2, \dots, l_s\}$ be a set of $s$ nonnegative integers. Suppose that $\mathcal{V}$ is a family of subspaces of $\mathbb{F}_{q}^{n}$ satisfying that $\dim(V_{i} \cap V_{j}) \in \mathcal{L}$ for any distinct subspaces $V_{i}$ and $V_{j}$ in $\mathcal{V}$ and
$\dim(V_{i}) \in \{k_{1}, k_{2}, \dots, k_{r}\}$ with $k_{i} > s - r$ for every $i$.
Then
\[|\mathcal{V}| \leq
\qbinom{n}{s}
+ \qbinom{n}{s-1} + \cdots + \qbinom{n}{s-r+1}.\]
\end{theorem}

We say that a family $\mathcal{V}$ of subspaces is {\em Sperner} (or antichain) if no subspace is contained in another subspace in $\mathcal{V}$, and $\mathcal{V}$ is
$k$-{\em Sperner} if all chains in $\mathcal{V}$ have length at most $k$.
The following well-known Sperner theorem for families of subspaces of $\mathbb{F}_{q}^{n}$ can be found in \cite{e2}, which is a vector space analogue of the classical Sperner theorem.

\begin{theorem}\label{thm1.10} (q-Analogue Sperner Theorem).  Assume that $\mathcal{V}$ is a Sperner family
of subspaces of $\mathbb{F}_{q}^{n}$. Then
\[|\mathcal{V}| \leq \qbinom{n}{\lfloor \frac{n}{2}\rfloor}.\]
\end{theorem}

We denote $\sum [n, k]$ to be the sum of the
$k$ largest $q$-binomial coefficients of order $n$, i.e., $\sum [n, k] = \sum_{i = 1}^{k} \left[n \atop \lfloor \frac{n - k}{2}\rfloor + i \right]$.
Let $\sum^{*} [n, k]$ be the collection of families consisting of the corresponding full levels,
i.e., if $n + k$ is odd, then $\sum^{*} [n, k]$ contains one family $\cup_{i = 1}^{k}\left[[n] \atop \lfloor \frac{n - k}{2}\rfloor + i \right]$
(where $\qbinom{[n]}{k}$ denotes the set of all $k$-dimensional subspaces of $\mathbb{F}_{q}^{n}$);
if $n + k$ is even, then $\sum^{*} [n, k]$ contains two families of the same size $\cup_{i = 0}^{k - 1}\left[[n] \atop \frac{n - k}{2} + i \right]$
and $\cup_{i = 1}^{k}\left[[n] \atop \frac{n - k}{2} + i \right]$.

The following $q$-analogue of Erd\H{o}s' theorem (Theorem \ref{thm1.5}) is implied by Theorem 2 in \cite{s2}, which generalizes Theorem \ref{thm1.10} to $k$-Sperner families.

\begin{theorem}\label{thm1.11} (Samotij \cite{s2}).  Assume that $\mathcal{V}$ is a $k$-Sperner family
of subspaces of $\mathbb{F}_{q}^{n}$. Then
\[|\mathcal{V}|\leq \sum [n, k].\]
Moreover, if $|\mathcal{V}| = \sum [n, k]$, then $\mathcal{V} \in \sum^{*} [n, k]$.
\end{theorem}

We say that a family $\mathcal{V}$ of subspaces is {\em intersecting} if $\dim(V_{i} \cap V_{j}) > 0$ for any pair of distinct $V_{i}, V_{j} \in \mathcal{V}$.
Here, our first main result is the following theorem on an intersecting $k$-Sperner family of subspaces of $\mathbb{F}_{q}^{n}$ which improves considerably the bound in Theorem \ref{thm1.11}, where $\mathcal{L}_{n, k}[R]$ denotes the set of all $k$-dimensional subspaces containing a given subspace $R$ in $\mathbb{F}_{q}^{n}$.

\begin{theorem}\label{thm1.12} Suppose that $\mathcal{V}$ is an intersecting $k$-Sperner family
of subspaces of $\mathbb{F}_{q}^{n}$ such that $\dim(V) \leq \lfloor \frac{n}{2} \rfloor$ for every  $V \in \mathcal{V}$. Then
\[|\mathcal{V}| \leq \sum_{j = \lfloor \frac{n}{2} \rfloor - k + 1}^{\lfloor \frac{n}{2} \rfloor}\qbinom{n-1}{j-1}\]
and equality holds if and only if there exists a $1$-dimensional subspace $R$ of $\mathbb{F}_{q}^{n}$ such that $\mathcal{V}_{j} = \{V \in \mathcal{V} \mid dim(V) = j\} = \mathcal{L}_{n, j}[R]$ for each $\lfloor \frac{n}{2} \rfloor - k + 1 \leq j \leq \lfloor \frac{n}{2} \rfloor$.
\end{theorem}

Note that for $\qbinom{n}{k}$ defined by (1), one has
\[\qbinom{n}{k} = q^{k}\qbinom{n-1}{k} + \qbinom{n-1}{k-1}.\]
The bound in Theorem \ref{thm1.12} is considerably smaller than that in Theorem \ref{thm1.11}.

We say that a property $P$ of families is {\em hereditary} if for any family $\mathcal{F}$ with property $P$,
every subfamily of $\mathcal{F}$ has property $P$. Clearly, $\mathcal{L}$-intersecting, Sperner, and $k$-Sperner are hereditary properties.

Our next main result establishes a general relationship between upper bounds for the sizes of $\mathcal{L}$-intersecting families of subsets of $[n]$ and  upper bounds for the sizes of $\mathcal{L}$-intersecting families of subspaces of $\mathbb{F}_{q}^{n}$.
Denote
\[f(n, c_{0}, \dots, c_{k}) = c_{k}{{n} \choose {k}}+ c_{k-1}{{n} \choose {k-1}}+ \cdots + c_{0}{{n} \choose {0}},\]
\[f_{q}(n, c_{0}, \dots, c_{k}) = c_{k}\qbinom{n}{k} +c_{k-1}\qbinom{n}{k-1} + \cdots + c_{0}\qbinom{n}{0},\]
where $0 \leq c_{i} \leq 1$ for each $0 \leq i \leq k-1$ and $0 < c_{k} \leq 1$.

\begin{theorem}\label{thm1.13}
Let $P$ be a hereditary property which is $\mathcal{L}$-intersecting or $T$-configuration forbidden, where $\mathcal{L}$ is a set of nonnegative integers and $T$ is a given configuration. Assume that $k$ is fixed and $n \geq 2k$.
If any family in $2^{[n]}$ satisfying property $P$ has at most
$f(n,  c_{0}, \dots, c_{k})$ members, then any family of subspaces of $\mathbb{F}_{q}^{n}$ with property $P$ has at most $f_{q}(n,  c_{0}, \dots, c_{k})$ members for $n$ sufficiently large.
\end{theorem}

By Theorem \ref{thm1.13}, one sees easily that for $n$ sufficiently large, Theorems \ref{thm1.2}--\ref{thm1.4} give rise to  Theorems \ref{thm1.7}--\ref{thm1.9}, respectively.

Let $0 < c \leq 1$ and
\[f^{*}(n, k, c) = c{{n} \choose {k}},\]
\[f^{*}_{q}(n, k, c) = c \qbinom{n}{k}.\]
For uniform families, we have the following similar result without requiring $n$ to be sufficiently large.

\begin{theorem}\label{thm1.14}
Let $P$ be a hereditary property which is $\mathcal{L}$-intersecting or $T$-configuration forbidden, where $\mathcal{L}$ is a set of nonnegative integers and $T$ is a given configuration. Assume that $k$ is fixed and $n \geq 2k$.
If any family in ${{[n]} \choose {k}}$ satisfying property $P$ has at most
$f^{*}(n, k, c)$ members, then any family of $k$-dimensional subspaces of $\mathbb{F}_{q}^{n}$ with property $P$ has at most $f^{*}_{q}(n, k, c)$ members.
\end{theorem}

A collection $M$ of pairwise disjoint sets is called a {\em matching}. The {\em matching number} $\nu(\mathcal{F})$ of
a family $\mathcal{F} \subseteq  {{[n]} \choose {k}}$ is the size of the largest matching that $\mathcal{F}$ contains.
Note that a family $\mathcal{F}$ is intersecting if $\nu(\mathcal{F}) = 1$. Erd\H{o}s \cite{e1} provided the following well-known conjecture about
the maximum size of a family $\mathcal{F} \subseteq  {{[n]} \choose {k}}$ with $\nu(\mathcal{F}) \leq s$.

\begin{conjecture}\label{conj1.15} (Erd\H{o}s Matching Conjecture \cite{e1}).
Let $\mathcal{F} \subseteq  {{[n]} \choose {k}}$ with $n \geq (s + 1)k$ and $\nu(\mathcal{F}) \leq s$.
Then
\[|\mathcal{F}| \leq \max \bigg\{{{(s+1)k -1} \choose {k}}, \mbox{ } {{n} \choose {k}} - {{n - s} \choose {k}} \bigg\}.\]
\end{conjecture}

In section 4, we provide applications of Theorems \ref{thm1.13} and \ref{thm1.14}.
In particular, we derive generalizations of the well known results about Erd\H{o}s-Chv\'atal simplex conjecture and Erd\H{o}s matching conjecture.

\section{Proof of Theorem 1.12}
First, we show that $\left[n \atop k \right]$ satisfies the following symmetric and unimodal properties similar to those of the binomial coefficients.

\begin{lemma}\label{lem2.1}
For $0 \leq k < l \leq \lfloor \frac{n}{2} \rfloor$, we have
\[{\rm (i)} \mbox{ } \qbinom{n}{k} = \qbinom{n}{n-k}, \hspace{120mm}\]
\[{\rm (ii)} \mbox{ } \qbinom{n}{k} < \qbinom{n}{l}.  \hspace{126mm}\]
\end{lemma}

\noindent {\bf Proof.}
For (i), it follows from the fact
\[\qbinom{n}{n-k} = \prod_{0 \leq i \leq n-k-1}\frac{q^{n-i}-1}{q^{n-k-i} - 1}\hspace{86mm}\]
\[= \frac{(q^{n}-1)(q^{n-1}-1)\cdots (q^{n-(n-k)+1}-1)}{(q^{n-k}-1)(q^{n-k-1}-1)\cdots (q - 1)} \hspace{42mm}\]
\[\hspace{10mm} = \frac{(q^{n}-1)(q^{n-1}-1)\cdots (q^{n-(n-k)+1}-1)}{(q^{n-k}-1)(q^{n-k-1}-1)\cdots (q - 1)} \cdot
\frac{(q^{n}-1)(q^{n-1}-1)\cdots (q-1)}{(q^{n}-1)(q^{n-1}-1)\cdots (q-1)}\]
\[ = \frac{(q^{n}-1)(q^{n-1}-1)\cdots (q^{n-k+1}-1)}{(q^{k}-1)(q^{k-1}-1)\cdots (q - 1)} = \qbinom{n}{k}. \hspace{35mm}\]

\noindent For (ii), it suffices to show that for $k + 1 \leq \lfloor \frac{n}{2} \rfloor$,
\[\qbinom{n}{k} < \qbinom{n}{k+1}.\]
Note that
\[\qbinom{n}{k+1} = \frac{q^{n - k} - 1}{q^{k+1} - 1}\qbinom{n}{k}.\]
For $k + 1 \leq \lfloor \frac{n}{2} \rfloor$, we have $n - k > k + 1$. It follows that
\[\hspace{62mm} \qbinom{n}{k} < \qbinom{n}{k+1}. \hspace{62mm} \hfill \Box   \]

\vspace{3mm}

The following LYM type inequality for intersecting Sperner families of subspaces of $\mathbb{F}_{q}^{n}$
is given in \cite{w1}. Recall that $\mathcal{L}_{n, k}[R]$ denotes the set of all $k$-dimensional subspaces containing a given subspace $R$ in $\mathbb{F}_{q}^{n}$.

\begin{theorem}\label{thm2.2} (Wang, \cite{w1}).
Suppose that $\mathcal{V}$ is an intersecting Sperner family of subspaces of $\mathbb{F}_{q}^{n}$.
Then
\[\sum_{k = 1}^{\lfloor \frac{n}{2} \rfloor}\frac{|\mathcal{V}_{k}|}{\qbinom{n-1}{k-1}} \leq 1,\]
where $\mathcal{V}_{k} = \{V \in \mathcal{V} \mid dim(V) = k\}$ for $1 \leq k \leq \lfloor \frac{n}{2} \rfloor$.
Equality holds if and only if $\mathcal{V} = \mathcal{L}_{n, k}[R]$ for some $1 \leq k \leq \lfloor \frac{n}{2} \rfloor$ and some $1$-dimensional subspace $R$ of $\mathbb{F}_{q}^{n}$.
\end{theorem}

\begin{lemma}\label{lem2.3}
If a family $\mathcal{V}$ of subspaces of $\mathbb{F}_{q}^{n}$ is intersecting k-Sperner, then $\mathcal{V}$ is the union of $k$ disjoint intersecting antichains.
\end{lemma}

\noindent {\bf Proof.}
Let $\mathcal{V} = \{V_{1}, V_{2}, \dots, V_{m}\}$ be a family of subspaces of $\mathbb{F}_{q}^{n}$.
Suppose that $\mathcal{V}$ is intersecting $k$-Sperner.
We now show that $\mathcal{V}$ is the union of $k$ intersecting antichains.
Define the required $k$ antichains recursively as follows:
Let $\mathcal{W}_{1}$ denote the family of all minimal subspaces in $\mathcal{V}$, and if $\mathcal{W}_{j}$ is defined for all $1 \leq j < i$,
then denote $\mathcal{W}_{i}$ to be the family of all minimal subspaces in $\mathcal{V} \setminus (\cup_{j = 1}^{i - 1}\mathcal{W}_{j})$. This process continues and stops at $i = h$.
Then $\mathcal{W}_{1}$, $\mathcal{W}_{2}$, \dots, $\mathcal{W}_{h}$ are disjoint intersecting antichains by definition. We claim that $h \leq k$. For every $V \in \mathcal{W}_{i}$,
there exists a $V' \in \mathcal{W}_{i - 1}$ such that $V' \subseteq V$.
Hence, if there exists a subspace in $\mathcal{W}_{k+1}$, then we would have a $(k + 1)$-chain in $\mathcal{V}$, a contradiction.
Thus, $h \leq k$ and $\mathcal{V}$ is the union of $k$ disjoint intersecting antichains.
$\hfill \Box$

\vspace{3mm}

As an easy consequence of Theorem \ref{thm2.2} and Lemma \ref{lem2.3}, one has the next extension of Theorem \ref{thm2.2} to intersecting $k$-Sperner families.

\begin{theorem}\label{thm2.4}
Suppose that $\mathcal{V}$ is an intersecting $k$-Sperner family of subspaces of $\mathbb{F}_{q}^{n}$.
For $0 \leq j \leq n$, let $\mathcal{V}_{j} = \{V \in \mathcal{V} \mid \dim(V) = j\}$. If $\dim(V) \leq \lfloor \frac{n}{2} \rfloor$ for each $V \in \mathcal{V}$,
then
\[\sum_{j = 1}^{\lfloor \frac{n}{2} \rfloor}\frac{|\mathcal{V}_{j}|}{\qbinom{n-1}{j-1}} \leq k.\]
\end{theorem}

\noindent {\bf Proof.}
Suppose that $\mathcal{V} = \{V_{1}, V_{2}, \dots, V_{m}\}$ is an intersecting $k$-Sperner family of $\mathbb{F}_{q}^{n}$. By Lemma \ref{lem2.3}, $\mathcal{V}$ is the union of $k$ intersecting antichains $\mathcal{W}_{1}$, $\mathcal{W}_{2}$, \dots, $\mathcal{W}_{k}$.
Applying Theorem \ref{thm2.2} on each of these antichains $\mathcal{W}_{i}$ and adding all the inequalities up, one obtains the desired result.
$\hfill \Box$

\vspace{3mm}

\begin{lemma}\label{lem2.5}
Suppose that for integers $1 \leq k \leq \lfloor \frac{n}{2} \rfloor$ and nonnegative real numbers $f_{1}, f_{2}, \dots, f_{\lfloor \frac{n}{2} \rfloor}$,
the following inequalities hold:
\[\sum_{j = 1}^{\lfloor \frac{n}{2} \rfloor}\frac{|f_{j}|}{\qbinom{n-1}{j-1}} \leq k,\]
\[f_{i} \leq \qbinom{n-1}{i-1} \mbox{ for each } 1 \leq i \leq \bigg\lfloor \frac{n}{2} \bigg\rfloor.\]
Then
\[\sum_{j = 1}^{\lfloor \frac{n}{2} \rfloor}f_{j} \leq \sum_{j = \lfloor \frac{n}{2} \rfloor - k + 1}^{\lfloor \frac{n}{2} \rfloor}\qbinom{n-1}{j-1}.\]
\end{lemma}

\noindent {\bf Proof.}
Consider the vector $f = (f_{1}, f_{2}, \dots, f_{\lfloor \frac{n}{2} \rfloor})$ which maximizes $\sum f_{i}$.
For $\qbinom{n-1}{j-1} < \qbinom{n-1}{i-1}$, the inequalities $f_{i} < \qbinom{n-1}{i-1}$ and $f_{j} > 0$ would lead to a bigger $\sum f_{i}$ by replacing $f_{i}$ and $f_{j}$ by $f_{i} + \varepsilon \qbinom{n-1}{i-1}$ and $f_{j} - \varepsilon \qbinom{n-1}{j-1}$, where $\varepsilon > 0$ is small enough, giving a contradiction.
It follows that we have either $f_{i} = \qbinom{n-1}{i-1}$ or  $f_{i} = 0$ for each $0 \leq i \leq \lfloor \frac{n}{2} \rfloor$.
Thus, by Lemma \ref{lem2.1}, the lemma follows.
$\hfill \Box$

\vspace{3mm}

\noindent {\bf Proof of Theorem 1.12.}
Suppose that $\mathcal{V}$ is an intersecting $k$-Sperner family of subspaces of $\mathbb{F}_{q}^{n}$ such that $\dim(V) \leq \lfloor \frac{n}{2} \rfloor$ for each $V \in \mathcal{V}$.
Let $\mathcal{V}_{j} = \{V \in \mathcal{V} \mid \dim(V) = j\}$ and
$f_{j} = |\mathcal{V}_{j}|$ for each $0 \leq j \leq \lfloor \frac{n}{2} \rfloor$. Since $\mathcal{V}$ is intersecting, $\mathcal{V}_{j}$ is intersecting and $j$-uniform for each $1 \leq j \leq \lfloor \frac{n}{2} \rfloor$. It follows from Theorem \ref{thm1.6} that $f_{j} = |\mathcal{V}_{j}| \leq \qbinom{n-1}{j-1}$ for every $1 \leq j \leq \lfloor \frac{n}{2} \rfloor$.
By Theorem \ref{thm2.4}, we have
\begin{equation}
\sum_{j = 1}^{\lfloor \frac{n}{2} \rfloor}\frac{|\mathcal{V}_{j}|}{\qbinom{n-1}{j-1}} \leq k.
\end{equation}
It follows from Lemma \ref{lem2.5} that
\[|\mathcal{V}|= \sum_{j = 1}^{\lfloor \frac{n}{2} \rfloor}f_{j} \leq \sum_{j = \lfloor \frac{n}{2} \rfloor - k + 1}^{\lfloor \frac{n}{2} \rfloor}\qbinom{n-1}{j-1}.\]
Suppose that the equality holds above. Then
it follows from the proof of Lemma \ref{lem2.5} that we must have $|\mathcal{V}_{i}| = f_{i} = \qbinom{n-1}{i-1}$ or $|\mathcal{V}_{i}| = f_{i} = 0$ for each $0 \leq i \leq \lfloor \frac{n}{2} \rfloor$. Thus, it follows from Theorem \ref{thm1.6} that $\mathcal{V}_{j} = \mathcal{L}_{n, j}[R_{j}]$ with $dim(R_{j}) = 1$ for each
$ \lfloor \frac{n}{2} \rfloor - k + 1 \leq j \leq \lfloor \frac{n}{2} \rfloor$.
Moreover, since $n \geq 2k+1$ and $\mathcal{V}$ is intersecting, we must have $R_{i} = R_{j} = R$ for all $ \lfloor \frac{n}{2} \rfloor - k + 1 \leq i, j \leq \lfloor \frac{n}{2} \rfloor$, where $R$ is a fixed $1$-dimensional subspace of $\mathbb{F}_{q}^{n}$. For otherwise,
if $R_{i} \neq R_{j}$ for some $\lfloor \frac{n}{2} \rfloor - k + 1 \leq i, j \leq \lfloor \frac{n}{2} \rfloor$ with $i \neq j$, then there exist
$W_{i} \in \mathcal{V}_{i}$ and $W_{j} \in \mathcal{V}_{j}$ such that $dim(W_{i} \cap W_{j}) = 0$, contradicting the assumption that $\mathcal{V}$ is intersecting.
$\hfill \Box$

\section{Proofs of Theorems 1.13 and 1.14}
We begin with the following lemma.

\begin{lemma}\label{lem3.1}
Let $k \geq 1$ be fixed and $n \geq 2k + 2$. For any $\varepsilon > 0$, there exists $n_{0}$ such that for $n \geq n_{0}$, we have
\[{\rm (i)} \mbox{ } \sum_{i = 0}^{k}{{n} \choose {i}} < \varepsilon{{n} \choose {k+1}}, \hspace{110mm}\]
\[{\rm (ii)} \mbox{ }  \sum_{i = 0}^{k}\qbinom{n}{i} < \varepsilon\qbinom{n}{k+1}.  \hspace{110mm}\]
\end{lemma}

\noindent {\bf Proof.} We first prove (i). Since
${{n} \choose {k+1}} = \frac{n - k}{k+1}{{n} \choose {k}}$, we have
\[\varepsilon{{n} \choose {k+1}} - (k + 1){{n} \choose {k}} = \bigg(\frac{\varepsilon(n - k)}{k+1} - (k + 1)\bigg){{n} \choose {k}}.\]
For $n > \frac{(k+1)^{2}}{\varepsilon} + k$, we have $\frac{\varepsilon(n - k)}{k+1} - (k + 1) > 0$.
Note that $\sum_{i = 0}^{k}{{n} \choose {i}} \leq (k+1){{n} \choose {k}}$. It follows that
\[\varepsilon{{n} \choose {k+1}} > (k + 1){{n} \choose {k}} \geq \sum_{i = 0}^{k}{{n} \choose {i}}.\]
For (ii), note that
\[\qbinom{n}{k+1} = \frac{q^{n - k} - 1}{q^{k+1} - 1}\qbinom{n}{k}.\]
Similar to (i), by applying Lemma \ref{lem2.1}, one has the desired inequality in (ii).
$\hfill \Box$

\vspace{3mm}

As a consequence of Lemma \ref{lem3.1}, one can derive the next fact easily. Recall that
\[f_{q}(n, c_{0}, \dots, c_{k}) = c_{k}\qbinom{n}{k} +c_{k-1}\qbinom{n}{k-1} + \cdots + c_{0}\qbinom{n}{0},\]
where $0 \leq c_{i} \leq 1$ for each $0 \leq i \leq k$.

\begin{lemma}\label{lem3.2} Let $k \geq 1$ be fixed and $n \geq 2k$.
For $n$ sufficiently large, if $f_{q}(n, b_{0}, \dots, b_{k}) > f_{q}(n, c_{0}, \dots, c_{k})$, then we have $b_{k} > c_{k}$ or there exists $j_{0} < k$ such that $b_{j} = c_{j}$ for all $j_{0} + 1 \leq j \leq k$ and $b_{j_{0}} > c_{j_{0}}$.
\end{lemma}

\noindent {\bf Proof.}
Suppose that $j_{0}$ is the largest integer between $0$ and $k$ such that $b_{j_{0}} \neq c_{j_{0}}$. Then we have either $b_{k} \neq c_{k}$ or $j_{0} < k$ and $b_{j} = c_{j}$ for all $j_{0} + 1 \leq j \leq k$.
By Lemma \ref{lem3.1} (ii), we must have either $b_{k} > c_{k}$ or $b_{j} = c_{j}$ for all $j_{0} + 1 \leq j \leq k$ and $b_{j_{0}} > c_{j_{0}}$ for $n$ sufficiently large.
$\hfill \Box$

\vspace{3mm}

Next, we introduce some notions and the following covering lemma by Gerbner \cite{g}.
Let $S$ be a set and $S = S_{0} \cup S_{1} \cup \cdots \cup S_{n}$ be a partition (In the discussions here, $S$ is either the set $2^{[n]}$ of all subsets of $[n]$ or the set of all subspaces of $\mathbb{F}_{q}^{n}$, and $S_{i}$ will be level $i$, i.e., $S_{i} = {{[n]} \choose {i}}$ or the set of all $i$-dimensional subspaces of $\mathbb{F}_{q}^{n}$, respectively).
Given a vector $t = (t_{0}, t_{1}, \dots, t_{n})$, we say a family $\Gamma$ of subsets of $S$ is a $t$-{\em covering family} of $S$ if for each $0 \leq i \leq n$, each member in $S_{i}$ is contained in exactly $t_{i}$ sets in the family $\Gamma$.

Given a family $\mathcal{F} \subseteq S$, let $f_{i} = |\mathcal{F} \cap S_{i}|$ and $(f_{0}, f_{1}, \dots, f_{n})$ is called the {\em profile vector} of $\mathcal{F}$.
For a weight $\overline{w} = (w_{0}, w_{1}, \dots, w_{n})$ and a family $\mathcal{F} \subseteq S$,
let $\overline{w}(\mathcal{F}) = \sum_{i = 0}^{n}w_{i}|\mathcal{F} \cap S_{i}|$.
Denote $\overline{w/t} = (w_{0}/t_{0}, w_{1}/t_{1}, \dots, w_{n}/t_{n})$. The following lemma is Lemma 2.1 in \cite{g}.

\begin{lemma}\label{lem3.3}
Let $P$ be a hereditary property of subsets (families) of $S$ and $\Gamma$ be a
$t$-covering family of $S$. Assume that there exists a real number $x$ such that for every $G \in \Gamma$, every subset $G'$ of $G$ with property $P$
has $\overline{w/t}(G') \leq x$. Then $\overline{w}(F) \leq |\Gamma|x$ for every $F \subseteq S$ with property $P$.
\end{lemma}

Let $S = \mathcal{L}_{n}(q)$ be the set of all subspaces of $\mathbb{F}_{q}^{n}$.
Denote $S = S_{0} \cup S_{1} \cup \cdots \cup S_{n}$ for a partition of $S$ such that $S_{i}$ is the set of all $i$-dimensional subspaces of $\mathbb{F}_{q}^{n}$.
We construct a $t$-covering family $\Gamma$ of $S$ such that every $\mathcal{G} \in \Gamma$ is a subfamily in $\mathcal{L}_{n}(q)$ isomorphic to Boolean lattice $\mathcal{B}_{n}$ (the set $2^{[n]}$ of all subsets of $[n]$ with the ordering being the containment relation)
as follows: Choose an arbitrary basis $B = \{v_{1}, \dots, v_{n}\}$ of $\mathbb{F}_{q}^{n}$
and let $\mathcal{G}_{B} = \{span(U) \mid U \subseteq B\}$, i.e., the family of all subspaces that
are generated by subsets of the vectors in $B$.
Obviously the function that maps $H \subseteq [n]$ to the
subspace $span\{v_{x} \mid x \in H\}$ keeps inclusion and intersection properties.
Let $\Gamma$ be the collection of the families $\mathcal{G}_{B}$ over all bases $B$, i.e.,
$\Gamma = \{\mathcal{G}_{B} \mid B \mbox{ is a basis of } \mathbb{F}_{q}^{n}\}$.
Denote
\begin{equation}
\alpha(q, n) = \frac{(q^{n} - 1)(q^{n} - q)(q^{n} - q^{2}) \cdots (q^{n} - q^{n-1})}{n!}.
\end{equation}
Let $t = (t_{0}, t_{1}, \dots, t_{n})$ be such that for $0 \leq i \leq n$,
\begin{equation}
t_{i} =  \frac{(q^{i} - 1)(q^{i} - q) \cdots (q^{i} - q^{i-1})(q^{n} - q^{i}) \cdots (q^{n} - q^{n-1})}{i!(n-i)!}.
\end{equation}

The next lemma follows from an easy counting.

\begin{lemma}\label{lem3.4}
Let $S = \mathcal{L}_{n}(q)$ be the set of all subspaces of $\mathbb{F}_{q}^{n}$
and $\Gamma = \{\mathcal{G}_{B} \mid B \mbox{ is a basis of } \mathbb{F}_{q}^{n}\}$. Then $|\Gamma| = \alpha(q, n)$ and $\Gamma$ is a $t$-covering of $S$ with
$t = (t_{0}, t_{1}, \dots, t_{n})$ given in (4).
\end{lemma}

\noindent {\bf Proof.} To show $|\Gamma| = \alpha(q, n)$, it suffices to show that $\mathbb{F}_{q}^{n}$ has $\alpha(q, n)$ different bases.
To obtain a basis $B = \{v_{1}, v_{2}, \dots, v_{n}\}$ for $\mathbb{F}_{q}^{n}$, one can choose $n$ nonzero linearly independent vectors in the order
$v_{1}, v_{2}, \dots, v_{n}$ so that $v_{j} \not \in span\{v_{1}, v_{2}, \dots, v_{j-1}\}$ as follows: There are $q^{n} - 1$ ways to choose $v_{1}$, and then there are $q^{n} - q$ ways to choose $v_{2}$, $q^{n} - q^{2}$ ways to choose $v_{3}$, \dots,
$q^{n} - q^{n-1}$ ways to choose $v_{n}$. Thus, there are $(q^{n} - 1)(q^{n} - q)(q^{n} - q^{2}) \cdots (q^{n} - q^{n-1})$ combined ways to form  $B = \{v_{1}, v_{2}, \dots, v_{n}\}$ (in order). Clearly, there are $n!$ ways (permutations on $B$) to form the same $B$ (without order). It follows that there are
\[\alpha(q, n) = \frac{(q^{n} - 1)(q^{n} - q)(q^{n} - q^{2}) \cdots (q^{n} - q^{n-1})}{n!}\]
different bases in $\mathbb{F}_{q}^{n}$.

For any $i$-dimensional subspace $V$ in $\mathbb{F}_{q}^{n}$, similar to the argument above, there are
\[\frac{(q^{i}-1)(q^{i}-q)\cdots (q^{i}-q^{i-1})}{i!}\]
different bases in $V$. For each basis $D$ in $V$, one can extend $D$ to a basis of $\mathbb{F}_{q}^{n}$ in
\[\frac{(q^{n}-q^{i})(q^{n}-q^{i+1})\cdots (q^{n}-q^{n-1})}{(n-i)!}\]
different ways. Therefore, $V$ is contained in
\[t_{i} =  \frac{(q^{i} - 1)(q^{i} - q) \cdots (q^{i} - q^{i-1})(q^{n} - q^{i}) \cdots (q^{n} - q^{n-1})}{i!(n-i)!}\]
different $\mathcal{G}_{B} \in \Gamma$.
It follows that $\Gamma$ is a $t$-covering of $S = \mathcal{L}_{n}(q)$.
$\hfill \Box$

\vspace{3mm}

The following proof for Theorem \ref{thm1.13} is motivated by the proof of Theorem 1.3 in \cite{g}.

\vspace{3mm}

\noindent {\bf Proof of Theorem 1.13.}
Assume that $\mathcal{V}$ is a family in $\mathcal{L}_{n}(q)$ satisfying property $P$ which is either $\mathcal{L}$-intersecting or $T$-configuration forbidden.
Let $\Gamma = \{\mathcal{G}_{B} \mid B \mbox{ is a basis of } \mathbb{F}_{q}^{n}\}$. By Lemma \ref{lem3.4},
$\Gamma$ is a $t$-covering of $S = \mathcal{L}_{n}(q)$.
Let $w_{i} = t_{i}$ for each $0 \leq i \leq n$. Then $\overline{w/t} = (1, 1, \dots, 1)$.
Since every $\mathcal{G} \in \Gamma$ is isomorphic to Boolean lattice $\mathcal{B}_{n}$, by
the assumption, the largest weight $\overline{w/t}(\mathcal{G}')$ among all subfamilies $\mathcal{G}' \subseteq \mathcal{G}$
satisfying $P$, i.e.,
the largest cardinality of such $\mathcal{G}' \subseteq \mathcal{G}$ is $f(n,  c_{0}, \dots, c_{k})$. By Lemma \ref{lem3.3}, we have
\begin{equation}
\overline{w}(\mathcal{V}) \leq |\Gamma|f(n,  c_{0}, \dots, c_{k}).
\end{equation}
To the contrary, suppose that
$|\mathcal{V}| > f_{q}(n, c_{0}, \dots, c_{k})$.
Then there exists $b = (b_{0}, b_{1}, \dots, b_{n})$ with $0 \leq b_{j} \leq 1$ for $0 \leq j \leq n$ such that
\[|\mathcal{V}| = g_{q}(n, b_{0}, \dots, b_{n}) = b_{n}\qbinom{n}{n} + b_{n-1}\qbinom{n}{n-1} + \cdots + b_{0}\qbinom{n}{0}\]
with
\[\bigg|\mathcal{V} \cap \qbinom{[n]}{j} \bigg| = b_{j}\qbinom{n}{j} \mbox{ for } 0 \leq j \leq n,\]
where $\qbinom{[n]}{j}$ denotes the set of all $j$-dimensional subspaces of $\mathbb{F}_{q}^{n}$.
It follows that
\[\overline{w}(\mathcal{V}) = \sum_{i=0}^{n}b_{i}w_{i}\qbinom{n}{i} \hspace{95mm}\]
\[= \sum_{i=0}^{n}b_{i}\frac{(q^{i} - 1)(q^{i} - q) \cdots (q^{i} - q^{i-1})(q^{n} - q^{i}) \cdots (q^{n} - q^{n-1})}{i!(n-i)!}\qbinom{n}{i}  \]
\[= \sum_{i=0}^{n}b_{i} \frac{(q^{i} - 1) \cdots (q^{i} - q^{i-1})(q^{n} - q^{i}) \cdots (q^{n} - q^{n-1})}{i!(n-i)!} \hspace{18mm} \]
\[\hspace{66mm} \cdot \frac{(q^{n} - 1)(q^{n} - q)\cdots (q^{n} - q^{n-1})}{(q^{i} - 1) \cdots (q^{i} - q^{i-1})(q^{n} - q^{i}) \cdots (q^{n} - q^{n-1})}\]
\[ = \sum_{i=0}^{n}b_{i}\frac{\alpha(n, q)n!}{i!(n-i)!} = \sum_{i=0}^{h}b_{i}|\Gamma|{{n} \choose {i}} \hspace{51mm}\]
\[= g(n, b_{0}, \dots, b_{n})|\Gamma|, \hspace{72mm} \]
where $g(n, b_{0}, \dots, b_{n}) = \sum_{i=0}^{n}b_{i}{{n} \choose {i}}$. Combining with (5), we obtain
\begin{equation}
g(n, b_{0}, \dots, b_{n}) \leq f(n,  c_{0}, \dots, c_{k}).
\end{equation}
Since ${{n} \choose {j}} = {{n} \choose {n-j}}$ and $\qbinom{n}{j} =  \qbinom{n}{n - j}$  for $0 \leq j \leq \lfloor\frac{n}{2}\rfloor$ by Lemma \ref{lem2.1},
we may combine symmetric $j$-th and $(n-j)$-th terms in $g_{q}(n, b_{0}, \dots, b_{n})$ and $g(n, b_{0}, \dots, b_{n}) = \sum_{i=0}^{n}b_{i}{{n} \choose {i}}$
for $0 \leq j \leq \lfloor\frac{n}{2}\rfloor$, and obtain
\[|\mathcal{V}| = g_{q}(n, b_{0}, \dots, b_{n}) = g_{q}(n, b'_{0}, \dots, b'_{\lfloor\frac{n}{2}\rfloor}) = b'_{\lfloor\frac{n}{2}\rfloor}\qbinom{n}{\lfloor\frac{n}{2}\rfloor} + b'_{\lfloor\frac{n}{2}\rfloor-1}\qbinom{n}{\lfloor\frac{n}{2}\rfloor-1} + \cdots + b'_{0}\qbinom{n}{0},\]
\[g(n, b_{0}, \dots, b_{n}) = g(n, b'_{0}, \dots, b'_{\lfloor\frac{n}{2}\rfloor}) = b'_{\lfloor\frac{n}{2}\rfloor}{{n} \choose {\lfloor\frac{n}{2}\rfloor}}
+ b'_{\lfloor\frac{n}{2}\rfloor - 1}{{n} \choose {\lfloor\frac{n}{2}\rfloor - 1}} + \dots + b'_{0}{{n} \choose {0}},\]
where $0 \leq b_{j} \leq b'_{j} = b_{j} + b_{n-j} \leq 2$ for each $0 \leq j \leq \lfloor\frac{n}{2}\rfloor$ with an exception that $b'_{\frac{n}{2}} = b_{\frac{n}{2}}$ when $n$ is even.
Assume that $h$ is the largest integer such that $h \leq \lfloor\frac{n}{2}\rfloor$  and $b'_{h} > 0$. Then we have
\[|\mathcal{V}| = g_{q}(n, b_{0}, \dots, b_{n}) = g_{q}(n, b'_{0}, \dots, b'_{h}) > f_{q}(n, c_{0}, \dots, c_{k}).\]
By Lemma \ref{lem3.1} (ii), we must have either $h > k$ or $h = k$ and
\[|\mathcal{V}| = g_{q}(n, b'_{0}, \dots, b'_{k}) > f_{q}(n, c_{0}, \dots, c_{k}).\]
We consider the following two cases:

\noindent {\bf Case 1.} $h = k$. We have $g_{q}(n, b'_{0}, \dots, b'_{k}) = f_{q}(n, b'_{0}, \dots, b'_{k})  > f_{q}(n, c_{0}, \dots, c_{k})$. It follows from Lemma \ref{lem3.2} that either $b'_{k} > c_{k}$ or there exists $j_{0} < k$ such that $b'_{j} = c_{j}$ for all $j_{0} + 1 \leq j \leq k$ and $b'_{j_{0}} > c_{j_{0}}$. Thus, by Lemma \ref{lem3.1} (i), we have that for $n$ sufficiently large,
\[g(n, b_{0}, \dots, b_{n}) = g(n, b'_{0}, \dots, b'_{k}) = \sum_{i=0}^{k}b'_{i}{{n} \choose {i}} = f(n, b'_{0}, \dots, b'_{k}) > f(n,  c_{0}, \dots, c_{k}),\]
contradicting (6).

\noindent {\bf Case 2.} $h > k$. Since $k+1 \leq h \leq \lfloor\frac{n}{2}\rfloor$, we have ${{n} \choose {h}} \geq {{n} \choose {k+1}}$. It follows from
Lemma \ref{lem3.1} (i) that
\[g(n,  b_{0}, \dots, b_{n}) = g(n, b'_{0}, \dots, b'_{h}) \geq b'_{h}{{n} \choose {k+1}} > f(n,  c_{0}, \dots, c_{k})\]
for $n$ sufficiently large, contradicting (6) again.

Therefore, we have
\[|\mathcal{V}| \leq f_{q}(n, c_{0}, \dots, c_{k}). \]
$\hfill \Box$

\vspace{3mm}

We remark that the inequality (6) in the proof of Theorem \ref{thm1.13} above can be obtained without using Lemma \ref{lem3.3} as follows:
By Lemma \ref{lem3.4}, $\Gamma = \{\mathcal{G}_{B} \mid B \mbox{ is a basis of } \mathbb{F}_{q}^{n}\}$
is a $t$-covering of $S = \mathcal{L}_{n}(q)$, i.e., every $i$-dimensional subspace in $S_{i} = \qbinom{[n]}{i}$ is contained in $t_{i}$ members in $\Gamma$.
Denote $\mathcal{V}_{j} = \mathcal{V} \cap  \qbinom{[n]}{j}$ for each $0 \leq j \leq n$. Let
\[|\mathcal{V}| = g_{q}(n, b_{0}, \dots, b_{n}) = b_{n}\qbinom{n}{n} + b_{n-1}\qbinom{n}{n-1} + \cdots + b_{0}\qbinom{n}{0}\]
with
\[|\mathcal{V}_{j}| = \bigg|\mathcal{V} \cap \qbinom{[n]}{j} \bigg| = b_{j}\qbinom{n}{j} \mbox{ for } 0 \leq j \leq n,\]
where $0 \leq b_{j} \leq 1$ for every $0 \leq j \leq n$.
Since $\mathcal{G}_{B}$ is isomorphic to Boolean lattice $\mathcal{B}_{n}$ for each basis $B$ of $\mathbb{F}_{q}^{n}$, we have
\[\sum_{\mathcal{G}_{B} \in \Gamma}|\mathcal{V}\cap \mathcal{G}_{B}| \leq |\Gamma|\cdot f(n, c_{0}, \dots, c_{k}).\]
On the other hand, since every $i$-dimensional subspace in $S_{i}$ is contained in $t_{i}$ members in $\Gamma$,
it follows that
\[\sum_{\mathcal{G}_{B} \in \Gamma}|\mathcal{V}\cap \mathcal{G}_{B}|
= \sum_{\mathcal{G}_{B} \in \Gamma}\sum_{i=0}^{n}|\mathcal{V}_{i}\cap \mathcal{G}_{B}|
= \sum_{i=0}^{n} \sum_{\mathcal{G}_{B} \in \Gamma}|\mathcal{V}_{i}\cap \mathcal{G}_{B}| \hspace{12mm}\]
\[\hspace{19mm} = \sum_{i=0}^{n} |\mathcal{V}_{i}|\cdot t_{i}
= \sum_{i=0}^{n} b_{i}\qbinom{n}{i}\cdot t_{i}
= \sum_{i=0}^{n} b_{i}|\Gamma|{{n} \choose {i}}, \]
where the last equality is due to the following double counting
\[\qbinom{n}{i}\cdot t_{i} = \bigg|\{(V, \mathcal{G}_{B}) \mid \dim(V) = i, \mathcal{G}_{B} \in \Gamma \mbox{ and } V \in \mathcal{G}_{B}\}\bigg|  = |\Gamma| {{n} \choose {i}}.\]
Combining the expressions above, we have
\[g(n, b_{0}, \dots, b_{n}) = \sum_{i=0}^{n}b_{i}{{n} \choose {i}} \leq f(n, c_{0}, \dots, c_{k}).\]

\vspace{3mm}

\noindent {\bf Proof of Theorem 1.14.}
Assume that $\mathcal{V}$ is a family of $k$-dimensional subspaces of $\mathbb{F}_{q}^{n}$ satisfying property $P$. Let
\[|\mathcal{V}| = g_{q}(n, k, a) = a\qbinom{n}{k} \mbox{ with } 0 < a \leq 1.\]
Similar to the proof of Theorem \ref{thm1.13} above, we have
\begin{equation}
g(n, k, a) = a{{n} \choose {k}} \leq f^{*}(n, k, c) = c{{n} \choose {k}}.
\end{equation}
To the contrary, suppose that
\[|\mathcal{V}| = g_{q}(n, k, a) = a\qbinom{n}{k} > f^{*}_{q}(n, k, c) = c \qbinom{n}{k}.\]
Then we have $a > c$ which implies
\[g(n, k, a) = a{{n} \choose {k}} > f^{*}(n, k, c) = c{{n} \choose {k}},\]
contradicting (7). Therefore, the theorem follows.
$\hfill \Box$

\section{Applications}
Clearly, by applying Theorem \ref{thm1.13}, one concludes that, for $n$ sufficiently large, Theorems \ref{thm1.7}--\ref{thm1.9} follow from
Theorems \ref{thm1.2}--\ref{thm1.4}, respectively.

Next, we provide some applications of Theorem \ref{thm1.14}.

Note that
\begin{equation}
{{n-1} \choose {k-1}} = \frac{k}{n}{{n} \choose {k}},
\end{equation}
\begin{equation}
\qbinom{n-1}{k-1} = \frac{q^{k} - 1}{q^{n} - 1} \qbinom{n}{k}.
\end{equation}

\subsection{Erd\H{o}s-Chv\'atal Simplex Conjecture}

A $d$-$simplex$ is defined to be a collection $A_{1}, \dots, A_{d+1}$ of subsets of size $k$ of $[n]$ such that the intersection of all of them is empty, but the intersection of any $d$ of them is non-empty. A $d$-$cluster$ is a collection $A_{1}, \dots, A_{d+1}$ of subsets of size $k$ of $[n]$ such that
$\cap_{i = 1}^{d+1} A_{i} = \emptyset$ and $|\cup_{i = 1}^{d+1} A_{i}| \leq 2k$. If $\{A_{1}, \dots, A_{d+1}\}$
is both a $d$-simplex and a $d$-cluster, we say that it is a $d$-simplex-cluster.

In 1971, Erd\H{o}s conjectured that a family $\mathcal{F}$ of $k$-subsets of $[n]$ with no $2$-simplex (also known as a triangle) satisfies
\[|\mathcal{F}| \leq {{n - 1} \choose {k-1}}.\]
In 1974, Chv\'atal \cite{c} generalized Erd\H{o}s' conjecture as follows.

\begin{conjecture}\label{conj4.1} (Erd\H{o}s-Chv\'atal Simplex Conjecture \cite{c}).
Suppose $k \geq d + 1 \geq 3$, $n \geq \frac{k(d + 1)}{d}$, and  $\mathcal{F}$ is a family of $k$-subsets of $[n]$ with no $d$-simplex. Then
\[|\mathcal{F}| \leq {{n - 1} \choose {k-1}} = {{n} \choose {k}} - {{n - 1} \choose {k}},\]
with equality only if $\mathcal{F}$ is a star.
\end{conjecture}

In 2005, Mubayi and Verstra\"ete \cite{mv} proved that the Erd\H{o}s-Chv\'atal Simplex Conjecture is true for $d = 2$, and
in 2021, Currier \cite{c1} proved that the conjecture holds when $d \geq 3$ and $n \geq 2k - d + 2$ as shown below.

\begin{theorem}\label{thm4.2} (Mubayi and Verstra\"ete \cite{mv}, Currier \cite{c1}).
Suppose $k \geq d + 1 \geq 3$, $n \geq 2k - d + 2$, and  $\mathcal{F}$ is a family of $k$-subsets of $[n]$ with no $d$-simplex-cluster. Then
\[|\mathcal{F}| \leq {{n - 1} \choose {k-1}} = {{n} \choose {k}} - {{n - 1} \choose {k}},\]
with equality if $\mathcal{F}$ is a star.
\end{theorem}

Let $\mathcal{A} = \{A_{1}, A_{2}, \dots, A_{d+1}\}$ be a family of $k$-dimensional subspaces of $\mathbb{F}_{q}^{n}$.
We say that $\mathcal{A}$ is a $d$-$simplex$ if $\dim(\cap_{i = 1}^{d+1}A_{i}) = 0$ but $\dim(\cap_{j = 1}^{d}A_{i_{j}}) > 0$ for any $d$ distinct members
$A_{i_{1}}, \dots, A_{i_{d}} \in \mathcal{A}$. A family $\mathcal{A} = \{A_{1}, A_{2}, \dots, A_{d+1}\}$ of subspaces of $\mathbb{F}^{n}_{q}$ is a $d$-$cluster$ if
$\dim(\cap_{i = 1}^{d+1} A_{i}) = 0$ and $\dim(span(\cup_{i = 1}^{d+1} A_{i})) \leq 2k$.  A family $\{A_{1}, \dots, A_{d+1}\}$ is a $d$-simplex-cluster if it
is both a $d$-simplex and a $d$-cluster.

Recall that
\[\qbinom{n}{k} = q^{k}\qbinom{n-1}{k} + \qbinom{n-1}{k-1} =  \qbinom{n-1}{k} + q^{n-k}\qbinom{n-1}{k-1}.\]
The following asymptotic $q$-analogue of Theorem \ref{thm4.2} follows easily from Theorems \ref{thm1.14} and \ref{thm4.2}.

\begin{theorem}\label{thm4.3}
Let $k \geq d + 1 \geq 3$. Suppose that $\mathcal{V}$ is a family of $k$-dimensional subspaces of $\mathbb{F}_{q}^{n}$ with no $d$-simplex-cluster.
Then for $n$ sufficiently large,
\[|\mathcal{V}| < q^{n-k}\qbinom{n-1}{k-1} = \qbinom{n}{k} - \qbinom{n-1}{k}.\]
\end{theorem}

\vspace{3mm}

\noindent {\bf Proof.} First, note that forbidding $d$-simplex-cluster is a hereditary property. Using the identity (8) and applying Theorems \ref{thm1.14} and \ref{thm4.2} with
\[f^{*}(n, k, \frac{k}{n}) = \frac{k}{n}{{n} \choose {k}} = {{n - 1} \choose {k-1}},\]
we conclude that
\begin{equation}
|\mathcal{V}| \leq f^{*}_{q}(n, k, \frac{k}{n}) = \frac{k}{n}\qbinom{n}{k}.
\end{equation}
Since $\lim_{n \longrightarrow \infty}\frac{k}{n}\frac{q^{n} - 1}{(q^{k} - 1)q^{n-k}} = 0$,
it follows from (9) and (10) that for $n$ sufficiently large,
\[|\mathcal{V}| \leq \frac{k}{n}\qbinom{n}{k} = \frac{k}{n} \cdot \frac{q^{n} - 1}{q^{k} - 1}\qbinom{n-1}{k-1} \hspace{34mm}\]
\[= \frac{k}{n}\frac{q^{n} - 1}{(q^{k} - 1)q^{n-k}}q^{n-k}\qbinom{n-1}{k-1} < q^{n-k}\qbinom{n-1}{k-1}.\]
$\hfill \Box$

Recently, Liu \cite{l1} proved the next result which confirms a conjecture by Mubayi and Verstra\"ete \cite{mv}.

\begin{theorem}\label{thm4.4} (Liu \cite{l1}).
Let $d \geq k \geq 4$ and $n$ be sufficiently large. Suppose that  $\mathcal{F}$ is a family of $k$-subsets of $[n]$ with no non-trivial intersecting subfamily of size $d+1$. Then
\[|\mathcal{F}| \leq {{n - 1} \choose {k-1}} = {{n} \choose {k}} - {{n - 1} \choose {k}},\]
with equality if $\mathcal{F}$ is a star.
\end{theorem}

Since forbidden a non-trivial intersecting subfamily of size $d+1$ is a hereditary property, the following $q$-analogue of Theorem \ref{thm4.4} follows from Theorems \ref{thm1.14} and \ref{thm4.4} and (8) and (9).

\begin{theorem}\label{thm4.5}
Let $d \geq k \geq 4$ and $n$ be sufficiently large. Suppose that $\mathcal{V}$ is a family of $k$-dimensional subspaces of $\mathbb{F}_{q}^{n}$ with no non-trivial intersecting subfamily of size $d+1$.
Then
\[|\mathcal{V}| \leq  q^{n-k}\qbinom{n-1}{k-1} = \qbinom{n}{k} - \qbinom{n-1}{k}.\]
\end{theorem}

\subsection{Erd\H{o}s Matching Conjecture}
Recall that the {\em matching number} $\nu(\mathcal{F})$ of
a family $\mathcal{F} \subseteq  {{[n]} \choose {k}}$ is the size of the largest matching that $\mathcal{F}$ contains
and recall Erd\H{o}s Matching Conjecture (or EMC for short) from Conjecture \ref{conj1.15}.
Clearly, the case $s = 1$ of the EMC is the classical Erd\H{o}s-Ko-Rado theorem.
Bollob\'as, Daykin and Erd\H{o}s \cite{bde} established the EMC for $n\geq 2k^{3}s$; Huang, Loh and Sudakov \cite{hls} proved the EMC for $n \geq 3k^{2}s$;
and Frankl \cite{f2} showed that the EMC holds for $n \geq (2s + 1)k - s$.

\begin{theorem}\label{thm4.6} (Bollob\'as, Daykin and Erd\H{o}s \cite{bde}, Frankl \cite{f2}, Huang, Loh and Sudakov \cite{hls}).
Let $\mathcal{F} \subseteq  {{[n]} \choose {k}}$ with $\nu(\mathcal{F}) \leq s$ and $n \geq (2s + 1)k - s$.
Then
\[|\mathcal{F}| \leq  {{n} \choose {k}} - {{n - s} \choose {k}}.\]
\end{theorem}

We say that two subspaces are disjoint if their intersection is the trivial subspace (i.e., $0$-dimensional subspace). A collection of pairwise disjoint subspaces is called a {\em matching}.
Similar to the matching number $\nu(\mathcal{F})$ of
$\mathcal{F} \subseteq  {{[n]} \choose {k}}$,
one can define the matching number $\nu_{q}(\mathcal{V})$ for a family $\mathcal{V}$ of $k$-dimensional subspaces of $\mathbb{F}_{q}^{n}$ to be the size of the largest matching that $\mathcal{V}$ contains.

Observe that a family $\mathcal{F}$ has matching number $\nu(\mathcal{F}) \leq s$ is a hereditary property. By applying Theorem \ref{thm1.14}, we obtain the next asymptotic $q$-analogue of Theorem \ref{thm4.6}.
Note that for $n \geq 2k + s$ and $s \geq 1$, one has
\[\qbinom{n}{k} - \qbinom{n - 1}{k} \leq \qbinom{n}{k} - \qbinom{n - s}{k}.\]

\begin{theorem}\label{thm4.7}
Suppose that $\mathcal{V}$ is a family of $k$-dimensional subspaces of $\mathbb{F}_{q}^{n}$ with $\nu_{q}(\mathcal{V}) \leq s$ and $q \geq 2$.
Then for $n$ sufficiently large,
\[|\mathcal{V}| \leq q^{n-k}\qbinom{n-1}{k-1} = \qbinom{n}{k} - \qbinom{n - 1}{k}.\]
\end{theorem}

\vspace{3mm}

\noindent {\bf Proof.}
For $n \geq (2k+1)s$, by the identity (8), we have
\[{{n} \choose {k}} - {{n-s} \choose {k}} \hspace{100mm}\]
\[\hspace{20mm} = \bigg({{n} \choose {k}} - {{n-1} \choose {k}}\bigg) + \bigg({{n-1} \choose {k}} - {{n-2} \choose {k}}\bigg)
+ \cdots + \bigg({{n-s+1} \choose {k}} - {{n-s} \choose {k}}\bigg)\]
\[ = {{n-1} \choose {k-1}} + {{n-2} \choose {k-1}} + \cdots + {{n-s} \choose {k-1}} \leq s{{n-1} \choose {k-1}} = \frac{sk}{n}{{n} \choose {k}}. \hspace{13mm}\]
Since a family $\mathcal{F}$ has matching number $\nu(\mathcal{F}) \leq s$ is a hereditary property, by applying Theorems \ref{thm1.14} and \ref{thm4.6} with
\[f^{*}(n, k, \frac{sk}{n}) = \frac{sk}{n}{{n} \choose {k}},\]
we obtain
\begin{equation}
|\mathcal{V}| \leq f^{*}_{q}(n, k, \frac{sk}{n}) = \frac{sk}{n}\qbinom{n}{k}.
\end{equation}
Since $\lim_{n \longrightarrow \infty}\frac{sk}{n}\frac{q^{n} - 1}{(q^{k} - 1)q^{n-k}} = 0$,
it follows from (9) and (11) that for $n$ sufficiently large,
\[|\mathcal{V}| \leq \frac{sk}{n}\qbinom{n}{k} = \frac{sk}{n} \cdot \frac{q^{n} - 1}{q^{k} - 1}\qbinom{n-1}{k-1} \hspace{62mm}\]
\[= \frac{sk}{n}\frac{q^{n} - 1}{(q^{k} - 1)q^{n-k}}q^{n-k}\qbinom{n-1}{k-1} < q^{n-k}\qbinom{n-1}{k-1} = \qbinom{n}{k} - \qbinom{n - 1}{k}.\]
$\hfill \Box$

\vspace{3mm}

Recently, Ihringer \cite{i} provided the following conjecture about the largest $s$-EM-family of $k$-spaces.

\begin{conjecture} \label{conj4.8} (Ihringer \cite{i})
For $n \geq 2k$, if $Y$ is a largest $s$-EM-family of $k$-spaces of $\mathbb{F}_{q}^{n}$, then $Y$ is the union of $s$ intersecting families or its complement.
\end{conjecture}

Ihringer \cite{i} proved that this conjecture holds if $16s \leq \min \{q^{\frac{n-k-l+2}{3}}, q^{\frac{n-k-r}{3}}, q^{\frac{n}{2}-k+1}\}$.

The next conjecture given by Aharoni and Howard \cite{ah} and Huang, Loh and Sudakov \cite{hls}, independently, extends the Erd\H{o}s Matching Conjecture to a rainbow (multipartite) setting.

\begin{conjecture}\label{conj4.9} (Rainbow Erd\H{o}s Matching Conjecture \cite{ah} and \cite{hls}).
Let $n \geq (s+1)k$ and let $\mathcal{F}_{1}, \mathcal{F}_{2}, \dots, \mathcal{F}_{s+1} \subseteq {{[n]} \choose {k}}$.
If there are no pairwise disjoint $F_{1}, \dots, F_{s+1}$, where $F_{i} \in \mathcal{F}_{i}$ for every $i = 1, \dots, s+1$, then
\[\min_{i \in [s+1]} |\mathcal{F}_{i}| \leq  \max \bigg\{{{(s+1)k -1} \choose {k}}, \mbox{ } {{n} \choose {k}} - {{n - s} \choose {k}} \bigg\}.\]
\end{conjecture}

We call this conjecture the REMC for short.
Huang, Loh and Sudakov \cite{hls} proved that the REMC is true for $n \geq 3k^{2}s$, and recently, Lu, Wang and Yu \cite{lwy} verified the REMC for
$n > 2sk$ and $s$ sufficiently large, and Keller and Lifshitz \cite{kl} proved the REMC for $n > Ck$, where $C$ depends only on $s$.

\begin{theorem}\label{thm4.10} (Huang, Loh and Sudakov \cite{hls}, Keller and Lifshitz \cite{kl}, and Lu, Wang and Yu \cite{lwy}).
Assume that $\mathcal{F}_{1}, \mathcal{F}_{2},$
$\dots, \mathcal{F}_{s+1} \subseteq {{[n]} \choose {k}}$.
If there are no pairwise disjoint $F_{1}, \dots, F_{s+1}$, where $F_{i} \in \mathcal{F}_{i}$ for every $i = 1, \dots, s+1$ and $n > Ck$ with $C$ depending on $s$, then
\[\min_{i \in [s+1]} |\mathcal{F}_{i}| \leq  \max \bigg\{{{(s+1)k -1} \choose {k}}, \mbox{ } {{n} \choose {k}} - {{n - s} \choose {k}} \bigg\}.\]
\end{theorem}

Note that forbidding rainbow matching of size $s+1$ is a hereditary property.
By applying Theorems \ref{thm1.14} and \ref{thm4.10} and using the identities (8) and (9), we obtain the next asymptotic $q$-analogue of Theorem \ref{thm4.10}
similar to the proof of Theorem \ref{thm4.7}.

\begin{theorem}\label{thm4.11}
Assume that $\mathcal{V}_{1}, \mathcal{V}_{2}, \dots, \mathcal{V}_{s+1}$ are families of $k$-dimensional subspaces of $\mathbb{F}_{q}^{n}$.
If there are no pairwise disjoint $V_{1}, \dots, V_{s+1}$, where $V_{i} \in \mathcal{V}_{i}$ for every $i = 1, \dots, s+1$, then
\[\min_{i \in [s+1]} |\mathcal{V}_{i}| \leq \qbinom{n}{k} - \qbinom{n - 1}{k}\]
for $n$ sufficiently large.
\end{theorem}

\section{Concluding Remarks}
In section 2, we derived an upper bound for an intersecting $k$-Sperner family of subspaces of $\mathbb{F}_{q}^{n}$ which gives
a $q$-analogue of Erd\H{o}s' $k$-Sperner Theorem. In section 3, we provided a general connection between upper bounds for the sizes of $\mathcal{L}$-intersecting or $T$-configuration forbidden families of subsets of $[n]$ and upper bounds for the sizes of $\mathcal{L}$-intersecting or $T$-configuration forbidden families of subspaces of $\mathbb{F}_{q}^{n}$. As applications of Theorems \ref{thm1.13} and \ref{thm1.14}, we provided some results on families of subspaces of $\mathbb{F}_{q}^{n}$ which follow asymptotically from the corresponding theorems on families of subsets of $[n]$, including generalizations of existing results about the well-known Erd\H{o}s matching conjecture and Erd\H{o}s-Chv\'atal simplex conjecture.

We end this section with the following open problems which are $q$-analogues of Erd\H{o}s-Chv\'atal Simplex Conjecture and Erd\H{o}s Matching Conjecture.

\begin{conjecture}\label{conj5.1} ($q$-Analogue Erd\H{o}s-Chv\'atal Simplex Conjecture).
Let $k \geq d + 1 \geq 3$ and $n \geq \frac{k(d + 1)}{d}$.
Suppose that $\mathcal{V}$ is a family of $k$-dimensional subspaces of $\mathbb{F}_{q}^{n}$ with no $d$-simplex.
Then
\[|\mathcal{V}| \leq q^{n-k}\qbinom{n-1}{k-1} = \qbinom{n}{k} - \qbinom{n-1}{k}.\]
\end{conjecture}

\begin{conjecture}\label{conj5.2} ($q$-Analogue Erd\H{o}s Matching Conjecture).
Suppose that $\mathcal{V}$ is a family of $k$-dimensional subspaces of $\mathbb{F}_{q}^{n}$ with $n \geq (s + 1)k$ and $\nu_{q}(\mathcal{V}) \leq s$.
Then
\[|\mathcal{V}| \leq \min \bigg\{\qbinom{n}{k} - \qbinom{n - 1}{k}, \mbox{ }  s\qbinom{n-1}{k-1}\bigg\}\]
\end{conjecture}

Theorems \ref{thm4.3} and \ref{thm4.7} show that the conjectures above are true when $n$ is sufficiently large.

\section*{Declaration of Competing Interest}
The authors declare that they have no conflicts of interest to this work.

\section*{Acknowledgement}
The research is supported by the National Natural Science Foundation of China (71973103, 11861019, 12271527, 12071484).
Guizhou Talent Development Project in Science and Technology (KY[2018]046), Natural Science Foundation of Guizhou ([2019]1047, [2020]1Z001,  [2021]5609).

\end{document}